\documentclass{article}
\usepackage{amsfonts,amsmath,amsthm,stmaryrd}
\usepackage{graphicx}
\usepackage{pstricks}

\newcommand{\rem}[1]{}

\theoremstyle{plain}
\newtheorem{lemma}{Lemma}
\newtheorem{theorem}[lemma]{Theorem}

\newtheorem{proposition}[lemma]{Proposition}
\newtheorem{definition}[lemma]{Definition}

\theoremstyle{remark}
\newtheorem{remark}{Remark}

\newcommand*  {\N} {{\mathbb N}}

\newcommand  {\R} {{\mathbb R}}
\newcommand*  {\Z} {{\mathbb Z}}

\newcommand*{\abs}[1]{\lvert #1 \rvert}

\def\del{\partial}

\title{Mathematics for $2d$ Interfaces}

\author{
Claude Bardos\footnote{ Universit\'e Denis Diderot  and Laboratory JLL
Universit\'e Pierre et Marie Curie, Paris, France
Pauli fellow; Wolfgang Pauli Institute Vienna
(claude.bardos@gmail.com)} and 
David Lannes\footnote{ DMA, Ecole Normale Supérieure et CNRS UMR 8553, 45, rue d'Ulm, 75005 Paris.
(David.Lannes@ens.fr )}}
\date{ }

\begin{document}
\maketitle

\begin{abstract}
We present here a survey of recent results concerning the mathematical analysis of instabilities of the interface between two incompressible, non viscous, fluids of constant density and vorticity concentrated on the interface. This configuration includes the so-called Kelvin-Helmholtz (the two densities are equal), Rayleigh-Taylor (two different, nonzero, densities) and the water waves (one of the densities is zero) problems. After a brief review of results concerning strong and weak solutions of the Euler equation, we derive interface equations (such as the Birkhoff-Rott equation) that describe the motion of the interface. A linear analysis allows us to exhibit the main features of these equations (such as ellipticity properties); the consequences for the full, non linear, equations are then described. In particular, the solutions of the Kelvin-Helmholtz and Rayleigh-Taylor problems are necessarily analytic if they are above a certain threshold of regularity (a consequence is the illposedness of the initial value problem in a non analytic framework). We also say a few words on the phenomena that may occur below this regularity threshold. Finally, special attention is given to the water waves problem, which is much more stable than the Kelvin-Helmholtz and Rayleigh-Taylor configurations. Most of the results presented here are in $2d$ (the interface has dimension one), but we give a brief description of similarities and differences in the $3d$ case.

\end{abstract}
\tableofcontents

\section{Introduction}

This contribution is a survey of recent results concerning  the mathematical analysis of instabilities in interfaces between two fluids of different densities and vorticity concentrated on the interface. In $2d$ this interface is a time dependent curve and in $3d$ a time dependent surface. This includes the two following extreme  cases: 
\begin{itemize}
\item
The density of one of these fluids is infinitely small with respect to the density of the other and therefore can be taken equal to $0$; this corresponds to the study of water waves  or droplets in vacuum.
\item The two fluids have the same density (in fact there is only one fluid) and the vorticity is concentrated on an interface. 
\end{itemize}

Experiments  and numerical simulation show that  the range of possible behaviors for the interface is very rich. In some cases (especially in the water waves case) the interface is
 ``stable" and its description has led to many studies in both the physical and mathematical literature. In other cases (this is the generic case when both densities are nonzero), the interface is ``unstable". Our purpose is to show how mathematical analysis can contribute to the understanding of such instabilities. To do so, one starts from the most abstract model which is the initial value (or Cauchy) problem for the incompressible Euler equation. It models a configuration where viscosity and surface tension are neglected and therefore  it is already too abstract to describe physical situations. As recalled below it admits solutions that would correspond to instantaneous creation or extinction of energy and as such would solve the energy crisis. 
However since in many cases the surface tension and viscosity are small, it turns out to be essential to study this limit case.
We will consider the evolution of a fluid in a domain $\Omega $ which can be either the whole space $\Omega= \R^d$ a sub-domain $\Omega \subset \R^d$ or a periodic box $\Omega=(\R/\Z)^d$. The natural quantities are the density $\rho$ the velocity $u(x,t)$ and the pressure $p(x,t)$. The equation of motion in $\Omega$ are then given by the Euler equations (or conservation of momentum),
\begin{equation}\label{Naples1}
\partial_t (\rho u)+\nabla\cdot(\rho u\otimes u)+\nabla p=\rho \vec {\mathbf g},
\end{equation} 
where the gravity $\vec {\mathbf g}\in \R^d$ is a constant vector, and the conservation of mass,
\begin{equation}\label{consmass}
\partial_t \rho+\nabla\cdot (\rho u)=0.
\end{equation}
The fluid is also assumed to be incompressible; together with the conservation of mass,
this yields the condition
 \begin{equation}\label{incomp}
\nabla\cdot u=  0,
\end{equation} 
which can be viewed as a constraint on the evolution of the fluid with the pressure as the Lagrange multiplier of this constraint. In the presence of boundary ($\Omega \not= \R^d$) one has for all times $t$ the following impermeability condition ,
\begin{equation}\label{condN}
\forall x\in \del\Omega \qquad u(x,t)\cdot \vec n =0\, \hbox{ with }\quad \vec  n \hbox{ the exterior normal on } \del\Omega, \hbox{boundary  of }\Omega.
\end{equation}
If  $\Omega$ is unbounded, it is also assumed that
$
u(t,\cdot) $ vanishes at infinity.

 Besides the time variable   $t$, two   coordinate systems appear naturally.
In   the {\it Eulerian system}, $u(x,t)$ describes the velocity field at the point  $x$ and at the time  $t$. In the {\it Lagrangian system}, we rather use 
the quantity $\widetilde{u}(x,t)$ that stands for the velocity of the fluid particle that was at the point $x$ at the initial time $t=0$. Equivalently,
$\widetilde{u}(x,t)=u(X(x,t),t)$, where $t\mapsto X(x,t)$ solves  the ordinary differential equation:
\begin{equation}
\dot X(x,t) =u(X(x,t),t),\qquad X(x,t=0)=x \label{lagrange0}\,.
\end{equation}
We consider solutions $u(x,t)$ which are continuous  away from an interface $\Sigma(t)$ which defines a partition of $\Omega=
\Omega_\pm(t)\cup\Sigma(t)$. We assume that on both sides of this interface the density of fluid $\rho=\rho_\pm$ is constant.  It is important to keep in mind that the velocity of the fluid is well defined on $\Omega_\pm(t)$ and therefore that the equation (\ref{lagrange0}), together with the incompressibility condition (\ref{incomp}), defines two families of volume preserving diffeomorphisms from $\Omega_\pm(t=0)$ onto $\Omega_\pm(t)$.
However the velocity has no intrinsic definition on $\Sigma(t)$ . Only the global evolution of this interface matters; as a consequence there is some freedom in the choice of the parametrization $r(t,\lambda) $ of the interface $\Sigma(t)$, and this choice will play an important r\^ole below.

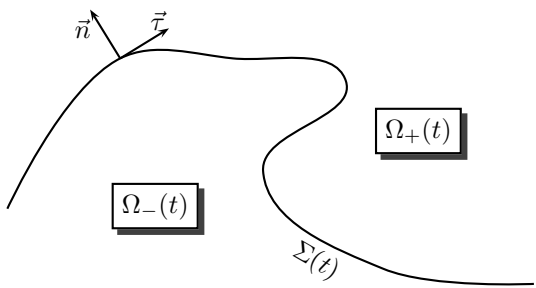
\begin{center}
\begin{pspicture}(0,1)(7,5)
\pscurve[showpoints=false](0,2)(1.5,4)(3,4)(4.5,3.7)(3.4,2.5)(5,1.2)(7,1)
\psline{->}(1.5,4)(1.1,4.65)
\rput(1,4.4){$\vec n$}
\psline{->}(1.5,4)(2.15,4.4)
\rput(2,4.5){$\vec\tau$}
\rput(5.5,3){\psshadowbox{$\Omega_+(t)$}}
\rput(2,2){\psshadowbox{$\Omega_-(t)$}}
\rput{330}(4.1,1.3){$\Sigma(t)$}
\end{pspicture}
\end{center}

Observe now that as long as the solution is smooth, the vorticity $\omega=\nabla\wedge u$ is solution of the system of equations:
\begin{equation}
\del_t \omega  + u  \cdot\nabla  \omega =\omega \cdot \nabla u  \,,\quad
\nabla\cdot u =0\,,\quad\nabla \wedge u  =\omega, \quad\hbox{ in } \Omega.
\end{equation}
In $2d$ with a planar flow, the vorticity can be written $\omega=(0,0,\partial_1 u_2-\partial_2 u_1)^T$ and the term $\omega\cdot\nabla u$ vanishes. Still denoting by
$\omega$ the \emph{scalar vorticity} $\omega=\partial_1 u_2-\partial_2 u_1$, 
the vorticity equation becomes
\begin{equation}\label{vort2d}
\del_t \omega  + u  \cdot\nabla  \omega =0.
\end{equation}
In any case, one has
\begin{equation}
\omega(x,0)=0\hbox{ for } x\notin \Sigma(0) \Rightarrow \omega(x,t)=0\hbox{ for } x\notin \Sigma(t),
\end{equation}
as long as the dynamics is well defined.
Several contributions have been devoted to the case where $\omega(0)$ is not $0$ away from $\Sigma(0)$ (cf. \cite{Lindblad,CoutandShkoller,ShatahZeng,ShatahZeng2,ZZ}). However since the evolution of the 
 interface is the dominant effect we will consider only dynamics with $0$ vorticity away from the interface. Therefore the only parameters are the Atwood number: 
  \begin{equation}
 a=  \frac {\rho_+-\rho_-}{\rho_++\rho_-}, \qquad -1 \le a \le 1,
 \end{equation}
 and the gravity $\vec{\mathbf g}$. The Rayleigh-Taylor problem corresponds to $0<\abs{a}<1\Leftrightarrow \rho_+\not=\rho_-\,, \rho_+\rho_-\not=0$, the Kelvin Helmholtz to $a=0\Leftrightarrow \rho_-=\rho_+$. Finally the situation where $a=\pm1\Leftrightarrow \rho_-\rho_+=0$ corresponds to a situation where one of the two densities vanishes and it describes  the dynamics of water waves. The Rayleigh-Taylor and the Kelvin-Helmholtz problems are highly unstable while the stability of the water waves problem is ensured if an extra condition ({\it called in the literature Taylor Hypothesis or Rayleigh-Taylor sign condition}, see \S \ref{rtcond}). The above names indicate that the study of these problems goes back to the end of the nineteenth century or the very beginning of the twentieth century. 
 
 The following facts should be kept in mind:
 \begin{itemize}
 \item The range of applications of the analysis is huge: From classical water waves (and application in oceanography)  to plasma physic in laser confinement and  astrophysics.
 \item The basic properties of the instabilities have already been established by the forerunners of the theory (which carries their names) by modal analysis of linearized problem near stationary solutions. A systematic use of the modal analysis can be found in the book of Chandrasekhar \cite{CHA}.
 \item The more recent contributions aim at a more systematic theory and allow for solutions that are not in a small neighborhood of stationary solutions.
\item We are interested here in the appearance of singularities and focus therefore on the most singular component of the solution (this is the central idea of symbolic analysis). However, the most singular parts are not necessarily the most important for the applications  (in the so called shallow water regime the terms neglected by the symbolic analysis are the most important ones, see \S \ref{shallow}).
\end{itemize}
For the sake of simplicity most of the presentation is devoted to  $2d$ problems. Similarities and differences which appear with the  $3d$ configurations are   considered in Section \ref{3d}.
 
 \section{Strong and weak solutions of the   $2d$ Euler equation}
 
When the density of both fluids is the same 
(say, $\rho+=\rho_-=1$, and thus Atwood number $a=0$), the
problem is reduced  to the ``genuine " incompressible Euler equations
(\ref{Naples1}), (\ref{incomp}). In this case, the gravity does not play any significant role (it can be put into the pressure term), and we thus take $\vec {\bf g}=0$ in (\ref{Naples1}):
\begin{equation}\label{Naples2}
\del_t u+\nabla\cdot( u\otimes u)+\nabla p=0,\qquad \nabla\cdot u=0.
\end{equation}
 By standard  methods \cite{Majda-Bertozzi} one can show the following local
existence theorem:
 \begin{theorem} For any initial data $u_0\in C^{1,\alpha }(\R^d)$ ($\alpha>0$, $d=2,3$)
with $\nabla\cdot u_0=0$ on $\R^d$,
there exists a finite time $T^*>0$ such that, with this initial data, one has  
a a unique classical solution to (\ref{Naples2})
on the time interval $[0,T^*]\,.$
\end{theorem}
 
 The above result is only local in time. $T^*$ depends on the size of the initial data and to the best of our knowledge this is the only ``general" result. Since we are concerned with solutions having their vorticity concentrated on a curve $\Sigma(t)$ (in $2d$, or a surface in $3d$), 
 the above theorem does not help much and weak solutions have to be considered
(for a precise meaning of what is meant by ``weak solution'', 
see for instance  \cite{Majda-Bertozzi}).
Using  the conservation of the vorticity along Lagrangian trajectories 
implied by (\ref{vort2d}), one gets the following classical results:
\begin{theorem} 1- (Yudovich \cite{YU}) With  an initial vorticity $\omega_0\in L^\infty(\R^2)\cap L^1(\R^2)$,  the Euler equations (\ref{Naples2})
 have a unique global weak solution which depends continuously on the data.\\
2-   If the initial vorticity is a bounded signed measure   \cite{DE}  or with a finite number of changes of sign  \cite{LNZ},  the Euler equations (\ref{Naples2}) have at least one weak solution $u(x,t)\in C_{\mathrm weak} (\R_t , L^2(\R^2))$.
\end{theorem}
However one should observe the big ``gap" between the hypothesis of the   point 1 and of the point 2 (this latter includes configurations allowing concentration of vorticity on an interface). Without the hypothesis  $\omega\in L^\infty(\R^2)$ no uniqueness or stability result is available. Furthermore {\it wild solutions}, for instance with space time compact support have been constructed by    Scheffer \cite{SC},  Shnirelman \cite{Sh} and finally  by  De Lellis and  Sz\'{e}kelyhidi  \cite{DeLellisandSzekelyhidi} (see also the nice review paper by Villani \cite{Villani}) who have obtained the following {\it instability, non uniqueness} theorem (note that the second point implies that the \emph{wild solutions} are not that exotic, since they are obtained as limits of strong, smooth solutions of the Euler equations):
\begin{theorem} 
 1-  For any bounded domain $O\subset\R^d_x\times \R_t$ ($d=2,3$) 
there exists a weak solution $u\in L^\infty(\R^d_x\times \R_t)$ to (\ref{Naples2}) which satisfies the relation  $$(x,t)\in O \Rightarrow |u(x,t)|=1\,, (x,t)\notin O \Rightarrow |u(x,t)|=0.$$
2- Moreover, $u$ is the limit in $L^2(\R^d_x\times \R_t)$ of a sequence of smooth, compactly supported, solutions $u_k$ ($k\in\N$) to the Euler equations with smooth (and vanishing as $k\to\infty$)
source terms $f_k$:
 $$\del_t u_k+  \nabla \cdot( u_k\otimes u_k)+\nabla p_k =f_k,
  \qquad \nabla\cdot u_k=0,\qquad f_k\to 0\quad \mbox{ in } \quad
{\mathcal D}'(\R^d).$$
  \end{theorem}

 \section{Weak solutions and interface equation}\label{weaksol0}
 In this section we derive a closed system of equations for the interface 
$\Sigma(t)$
 when  the vorticity is a   density concentrated on this curve,
assumed that it satisfies the following smoothness assumption:
 \begin{definition}
1- A rectifiable curve $\Sigma$ parametrized by its arc length, $\Sigma=\{r(s)\}$, is \emph{chord-arc} 
 if there exists two constants $0<m<M<\infty$ such that one has:
$$
m |s_1-s_2| \le | r(s_1)-r(s_2)| \le M|s_1-s_2|.
$$
A local version (near any point of $\Sigma$) of this definition is obtained with minor modifications.\\
2- We say that a curve  $\Sigma\subset \Omega$ is \emph{regular} if it is bounded away from $\del \Omega$, rectifiable and  chord-arc .
\end{definition}
 
 \subsection{Biot-Savart formula and Hilbert Transform}\label{bsht}

From the incompressibility relation $\nabla\cdot u=0$  one infers the
existence of a \emph{stream function} $\psi$ such that
\begin{equation}\label{BS1}
u(x,t)= \nabla^\perp  \psi (x,t)=\big(-\partial_y \psi(x,t),\partial_x\psi(x,t)\big)^T.
\end{equation} 
Classical results of potential theory show that $u(x,t)$ converges to two
values $u_\pm(x,t)$ on the two sides of the interface. The incompressibility 
condition implies the continuity across
the interface of the normal\footnote{It is of course assumed that the interface $\Sigma(t)$ is an orientable curve. We choose the orientation such that the normal vector $\vec n$ points from $\Omega_-(t)$ to $\Omega_+(t)$. The orientation of the tangent vector $\vec \tau$ is chosen such that $(\vec \tau,\vec n)$ is a direct orthonormal basis of $\R^2$} component of the velocity, but not of the tangential ones.
Denoting by $u^\tau_\pm$ the trace of the tangential component of the velocity on both sides of the interface, and by $r(\cdot,t)$ a parametrization\footnote{Changing $\lambda$ into $-\lambda$ if necessary, we always assume that $\partial_\lambda r$ is positively colinear to $\tau$, that is, $\vec \tau=\frac{1}{\vert\partial_\lambda r\vert}\partial_\lambda r$.} of the interface, $\Sigma(t)=\{r(\lambda,t)\}$, standard distribution calculus shows that the
scalar vorticity $\omega$ is a density measure supported on $\Sigma(t)$,
$$
\forall \varphi\in C(\Omega),\qquad\langle \varphi,\omega\rangle =
\int \varphi(r(\lambda,t))\widetilde{\omega}(\lambda,t)d\lambda,
\quad \mbox{ with }\quad
\widetilde{\omega}(\lambda,t)=-\abs{\partial_\lambda r}(u_+^\tau-u_-^\tau)_{\vert_{(r(\lambda,t),t)}};
$$
the quantity $\widetilde{\omega}$ is called the \emph{vorticity density} 
associated to the parametrization $r(\lambda,t)$. The vorticity density associated to the arc length parametrization is called the \emph{vortex strength}, and will be denoted $\gamma(s,t)$ throughout this text,
\begin{equation}\label{vortexstrength}
\gamma(s(\lambda,t),t)=\frac{1}{\abs{\partial_\lambda r(\lambda,t))}}\widetilde{\omega}(\lambda,t),
\quad \mbox{ with }\quad s(\lambda,t)=\int_0^\lambda \abs{\partial_\lambda r(\lambda',t)}d\lambda'.
\end{equation}
 In
terms of $\psi$, one has
\begin{equation}\label{BS2}
\Delta\psi=\omega=\widetilde\omega\otimes\delta_{\Sigma(t)}.
\end{equation}
From (\ref{BS1}), (\ref{BS2}) and the fact that $u$ vanishes at infinity if
$\Omega$ is unbounded, one can express $u(x,t)$ for all $x\notin \Sigma(t)$ in terms of the vorticity density $\widetilde{\omega}$
by a \emph{Biot-Savart} type formula. Denoting by $G(x,y)$ the Green function of the Laplacian in $\Omega$ with Dirichlet boundary condition, one has
\begin{eqnarray}
\forall x\in \Omega(t) \backslash  \Sigma (t), \qquad u(x,t)&=& \int_{\Omega}  \nabla^\perp_xG(x, y)\omega(t,y)d y \nonumber\\
&=& \int \nabla^\bot_x G(x,r(\lambda,t))\widetilde{\omega}(\lambda,t)
d \lambda\label{biotsavart1}\,.\end{eqnarray}
  For  $ x\in \Omega \backslash  \Sigma (t)$ 
the function  $u(t,x)$ is very smooth (indeed analytic). Though this is no more true for   $x\in \Sigma(t)$,  the formula (\ref{biotsavart1}) can be extended to the case $x\in \Sigma(t)$  provided the left hand side is replaced by 
$$
v= \frac{u_++u_-}2,
$$ 
and that on the right hand side the integral is defined as a principal value:
\begin{equation}
v(x,t)=\mathrm{p.v.} \int \nabla^\bot G(x,r(\lambda,t))\widetilde \omega(\lambda,t)d\lambda.\label{Biot-Savard}
\end{equation}
                                                                                       
In particular for $\Omega=\R^2$  one recovers the standard Biot-Savart formula:
\begin{equation}
v(x,t)=
\frac1{2\pi}R_{\frac{\pi}{2}} \mathrm{p.v.} \int \frac {x-r(\lambda,t)}
{|x-r(\lambda,t)|^2} \widetilde \omega(\lambda,t)  d\lambda,
 \label{biotsavard2}
\end{equation}
with $R_{\frac{\pi}{2}}$ denoting the rotation of angle $\frac\pi 2\,.
$
Furthermore, due to the ellipticity of the Laplacian and also for $\Omega \not=\R^2$, a local version of the above  formula remains valid modulo a more regular term 
\begin{eqnarray}
v(x,t)=\frac1{2\pi}R_{\frac{\pi}{2}} \mathrm{p.v.} \int_{\{r(\lambda,t)\in {\mathcal U}\}} \frac {x-r(\lambda,t)}
{|x-r(\lambda,t)|^2} \widetilde \omega(\lambda,t) d\lambda\label{biotsavard3} +R_{\mathcal U}(\omega);
\end{eqnarray}
in (\ref{biotsavard3}), $x\mapsto R_{\mathcal U}(\omega)(x) $ is an analytic function and $\omega\mapsto  R_{\mathcal U}(\omega)$ is a ``regularizing operator " which depends on $\Omega$ and the open set $\mathcal U$.

\begin{remark} \label{compacite}
Denoting by $v^\tau$ and $v^n$ the tangential and normal components of $v$,
one checks that:
\begin{eqnarray}
&&v^\tau(x,t)=-\frac1{2\pi} \mathrm{p.v.} \int_{\{r(\lambda,t)\in {\mathcal U}\}}
\vec n\cdot \frac {x-r(\lambda,t)}
{|x-r(\lambda,t)|^2} \widetilde \omega(\lambda,t)d\lambda\label{biotsavard4} \,
+R^\tau_{\mathcal U}(\omega)\\
&&v^n(x,t)=\frac1{2\pi} \mathrm{p.v.} \int_{\{r(\lambda,t)\in {\mathcal U}\}} 
\vec\tau\cdot\frac {x-r(\lambda,t)}
{|x-r(\lambda,t)|^2} \widetilde \omega(\lambda,t) d\lambda +R^n_{\mathcal U}(\omega)\label {biotsavard5}.
\end{eqnarray}
By a Taylor expansion of $r(\cdot,t)$ in the neighborhood of $x$, one checks
that  $(x-r(\lambda,t))\cdot \vec n=O(|x-r(\lambda,t)|^2)$ while one only has
 $(x-r(\lambda,t))\cdot \tau =O(|x-r(\lambda,t)|)$. It follows that
 $v^\tau$ is ``one degree " more regular that $v^n$.
\end{remark}
One of the main issues in the analysis of the Rayleigh-Taylor (and Kelvin-Helmholtz)  problem is the ellipticity of the linearized operator. This property is local by its definition and its effects. Therefore it is convenient to write, with a Galilean transform,  the equation  (\ref{biotsavard3}) in the neighborhood of a point $r(\lambda,t)=(0,0))$, with $\vec \tau$, $\vec n$ given by the coordinate axis, and where  $\Sigma(t)$ is a graph  $(x_1,\epsilon \sigma(x_1,t))$ which is a
 small perturbation of the $x_1$ axis in a neighborhood $\abs{x_1}<\delta$ of the origin. In this setting (\ref{biotsavard3}) becomes 
\begin{eqnarray}
&& v_1(x,t)=-\frac\epsilon {2\pi} \mathrm{p.v.}\int_{\abs{x_1}< \delta} \frac{\sigma(x_1,t)-\sigma(x'_1,t)}{(x_1-x_1')^2+\epsilon^2(\sigma(x_1,t)-\sigma(x'_1,t))^2}\tilde \omega(x_1',t)d  x_1',\label{bs1}\\
&& v_2(x,t)=\frac1{2\pi} \mathrm{p.v.}\int_{\abs{x_1}<\delta} \frac{x_1-x_1'}{(x_1-x_1')^2+\epsilon^2 (\sigma(x_1,t)-\sigma(x_1',t))^2}\tilde \omega(x_1',t)
d x_1',\label{bs2}
\end{eqnarray}
 where we omitted the more regular terms coming from $R_{\mathcal U}(\omega)$. We thus see that two operators will play an important role in the analysis of the linearized equation (see \S \ref{linan} below): the {\it Hilbert transform} $H$  and $|D|$, defined respectively as
\begin{eqnarray}
Hf(x_1)&=&\frac{1}{\pi} {\mathrm p.v.} \int\frac{1}{x_1-x_1'} f(x_1')dx_1'= F^{-1}(-i {\mathrm {sgn} }(\xi) \hat f (\xi))\label{ht1}\,,\\
|D|f(x_1)&=&\frac1{ \pi} {\mathrm p.v.} \int \frac{f(x_1)-f(x_1')}{(x_1-x_1')^2} = \del_{x_1}  (Hf(x_1))=  F^{-1}(|\xi|) \hat f (\xi))\,.\label{ht2}
\end{eqnarray}

\subsection{Jump relations and equations on the interface}\label{jump}

With minimal regularity hypothesis on the interface and on the vorticity density one can deduce ``jump relations". We have already seen that
 incompressibility relation $\nabla\cdot u=0$
and the identity $\nabla\wedge u=\widetilde\omega\otimes \delta_{\Sigma(t)}$
lead to 
\begin{equation}\label{firstjump}
[u^n]=0\quad\mbox{ and }\quad [u^\tau]=\frac{-1}{\abs{\partial_\lambda r}}\widetilde\omega,
\end{equation}
where we use the standard notation $[u]=u_+-u_-$.\\
From the conservation equations of mass (\ref{consmass}) and momentum (\ref{Naples1}),
written in the sense of distributions,  we deduce further two equations for the evolution of the interface.
If $r(\cdot,t)$ is a parametrization of the interface, 
$\Sigma(t)=\{r(\lambda,t)\}$, one gets (cf. \cite{SS}, and \cite{CCG} for a slightly different, but equivalent formulation):
 \begin{eqnarray}
(\del_tr -v)\cdot \vec n &=&0\label{tangent}\\
\del_t (\frac{1}{2}\widetilde\omega-a\abs{\partial_\lambda r}v^\tau)
&+&\partial_\lambda\Big(\frac{1}{\abs{\partial_\lambda r}}(\frac{1}{2}\widetilde\omega-a\abs{\partial_\lambda r} v^\tau)(v-\del_t r)\cdot\vec\tau\Big)\nonumber\\
&-&a\partial_\lambda\Big(\frac{1}{8}\frac{\widetilde{\omega}^2}{\abs{\partial_\lambda r}^2}-\frac{1}{2}\abs{v}^2\Big)=-a\abs{\partial_\lambda r}\vec{\bf g}\cdot\vec\tau.\label{vortex}
\end{eqnarray} 
Observe that the above system does not fully determines the evolution of the vector  $r(\lambda,t)$. The equation (\ref{tangent}) concerns only the normal component of its time derivative. This correspond to the fact already noticed that it is only the global evolution of the interface and not of its parametrization that matters. 
Different choices of the parametrization lead to the following formulas:
\begin{itemize}
\item
When $a=0$  (this corresponds to $\rho_-=\rho_+$ i.e. to the Kelvin-Helmholtz problem),  it is convenient to choose a Lagrangian parametrization, 
\begin{equation}
r_t(\lambda,t)= v(r(\lambda,t),t)=\frac{ u_++u_-}2(r(\lambda,t));
\end{equation}
the equations (\ref{tangent}) and (\ref{vortex}) then  reduce to
 \begin{equation} 
\del_t \tilde \omega =0
\quad \mbox{ and thus }\quad
 \tilde\omega(\lambda,t)= \tilde\omega(\lambda,0)\label{prebr1}.
\end{equation}
Using (\ref{Biot-Savard}), we thus get
\begin{equation}
\del_tr(\lambda,t)=\mathrm{p.v.}  \int \nabla^\bot_x G\big(r(\lambda,t), r(\lambda',t)\big)\tilde \omega(\lambda',0)d\lambda'
\label{prebr2}
\end{equation}
or, when $\Omega =\R^2$,
\begin{equation}
\del_t r(\lambda,t)= \frac1{2\pi}R_{\frac{\pi}{2}} \mathrm{p.v.} \int \frac {r(\lambda,t)-r(\lambda',t)}
{|r(\lambda,t)-r(\lambda',t)|^2} \tilde \omega(\lambda',0) d\lambda'\,.
\label{prebs0}
\end{equation}
It is therefore natural to reparametrize the interface by the (time independent) arc length \emph{of its initial position}, $\Sigma(t)=\{r(s_0^{-1}(s),t)\}$,
where $s_0(\lambda)=\int_0^\lambda \partial_\lambda r(\lambda',0)d\lambda'$.
 When $\Omega=\R^2$, one
can further identify the fluid domain with the complex plane, i.e. $z(s,t)\sim r(s_0^{-1}(s),t)$ to obtain:
 \begin{equation}
\del_t \overline z(s,t) = \frac{1}{2\pi i }\mathrm{p.v.}\int \frac{\gamma(s',0)}{z(s,t)-z(s',t) }ds',
\end{equation}
where $\gamma(s,0)$ is the vortex strength (\ref{vortexstrength})  at $t=0$.\\
If moreover $\gamma(s,0)$ has \emph{a distinguished sign} (which is the case when the
initial vorticity is a Radon measure), then we can use the mapping
$$
\alpha_0(s) =\int_0^s \gamma(s',0)ds'
$$
to perform another (time independent) reparametrization of the interface, namely, $\Sigma(t)=\{z(\alpha_0^{-1}(\alpha),t)\}$.
One thus obtains the \emph{Birkhoff-Rott equation} for the variable 
${\bf z}(\alpha,t)=z(\alpha_0^{-1}(\alpha),t)$,
 \begin{equation}\label{BR}
  \del_t \overline {\bf z}(\alpha,t) = \frac{1}{2\pi i }\mathrm{p.v.}\int \frac {d\alpha'} {{\bf z}(\alpha,t)-{\bf z}(\alpha',t) }. 
\end{equation}
The Birkhoff-Rott equation corresponds to a parametrization by the circulation; indeed, one
readily checks that $|\partial_\alpha{\bf z}(\alpha_0(s),t) |=\frac{1}{\gamma(s,0)}$. 
\item When $a=-1$ (this corresponds to $\rho_+=0$, and thus to the water waves problem), 
several choices are possible. One can for instance, as in \cite{SS}, choose a Lagrangian parametrization 
with the velocity given by the trace of the velocity field in $\Omega_-(r)$ at the surface,
\begin{equation}\label{paramWW}
r_t(\lambda,t)=u_-(r(\lambda,t),t)=u_-^n(r(\lambda,t),t)\vec n+u_-^\tau(r(\lambda,t),t)\vec \tau.
\end{equation}
Another choice, adopted in \cite{AM2}, consists in parametrizing the surface by its 
(renormalized) arclength. It is also possible to use the particular case $a=-1$ of the following point when the surface is a graph.
\item When the interface is a graph (we have seen that this hypothesis can always be used for ``local analysis"), we choose the canonical parametrization by $x_1$, $\Sigma(t)=\{(x_1,\sigma(x_1,t)\}$ and   the full system (including the Biot-Savart equation) becomes:
\begin{eqnarray}
&&\del_t \sigma + v_1\del_{x_1} \sigma = v_2\label{nlgtangent}\,\\
&& \del_t\big(\frac{1}{2}\widetilde\omega - a(v_1+v_2\del_{x_1}\sigma )\big)\nonumber \\ 
&&\qquad+\del_{x_1}\Big\{ v_1\big( \frac{1}{2}{\widetilde \omega} -a (v_1+v_2\del_{x_1}\sigma ) \big)
-a\big( \frac{\tilde \omega ^2}{8(1+\del_{x_1}\sigma^2)}- \frac{|v|^2}{2} -g_1 x_1-g_2 \sigma\big) \Big\}
  = 0 \label{nlgvortex}\,,\\
&& v_1(x_1,t)=-\frac 1 {2\pi}{\mathrm{p.v.}}\int \frac{\sigma(x_1,t)-\sigma(x_1',t)}{(x_1-x_1')^2+ (\sigma(x_1,t)-\sigma(x'_1,t))^2}\tilde \omega(x'_1,t)dx_1'\label{nlbs1}\,,\\
&& v_2(x_1,t)=\frac1{2\pi}{\mathrm{p.v.}}\int \frac{x_1-x_1'}{(x_1-x_1')^2+(\sigma(x_1,t)-\sigma(x_1',t))^2}\tilde \omega(x'_1,t)dx_1'\,.\label{nlbs2}
\end{eqnarray} 
\end{itemize}
\section{Linearization of the interface problem}\label{linan}
To analyze the local properties of the above system (\ref{nlgtangent})-(\ref{nlbs2}), one  introduces its  linearized version 
near a stationary solution $\underline{\sigma}(x,t)=0$, 
$\underline{u}_+(x,0)=(-u_0,0)$, and
 $\underline{u}_-(x,0) =(u_0,0)$ (and thus $\underline{\widetilde\omega}=2u_0$). Looking for solutions of the form $\underline{\sigma}+\epsilon \sigma$, $\underline{\widetilde\omega}+\epsilon\widetilde\omega$, we are led to the following equations:
\begin{eqnarray}
&&\del_t \sigma- v_2=0\,, \label{lt3} \\
&&\del_t ( \frac{1}{2}\widetilde\omega- av_1)+ \del_{x_1}\big(  u_0 v_1- \frac{1}{2}a u_0 \widetilde\omega   +ag_2\sigma\big)=0 \label{lv3}\,,\\
&& \epsilon v_1(x_1,t)=-\epsilon \frac1{2\pi}\mathrm{p.v.} \int \frac{\sigma(x_1,t)-\sigma(x'_1,t)}{(x_1-x_1')^2} (2u_0+\epsilon  \widetilde\omega(x_1',t))dx_1'\label{lbs3}=-\epsilon  u_0 |D|\sigma +O(\epsilon^2)\,,\\
&& \epsilon v_2(x_1,t)=\frac1{2\pi} \mathrm{p.v.}\int \frac{x_1-x_1'}{(x_1-x_1')^2 } (2u_0+\epsilon  \omega(x_1',t))dx_1'
=\frac{1}{2}\epsilon H\widetilde\omega +O(\epsilon^2),
\label{lbs4}\\
\end{eqnarray}
where the operators $\abs{D}$ and $H$ are as defined in (\ref{ht1}) and (\ref{ht2}) (and where the identity $\abs{D}1=H1=0$ has been used to derive the second part of (\ref{lbs3}) and (\ref{lbs4})).

At leading order in $\epsilon$, we thus obtain the two following equations:
\begin{eqnarray}
&&\del_t \sigma- \frac{1}{2}H\widetilde\omega=0 \label{lt2.3} \,,\\
&&\del_t ( \frac{1}{2}\widetilde\omega+ au_0|D|\sigma)+ \del_{x_1}\big(- u_0^2\abs{D}\sigma- \frac{1}{2}au_0\widetilde\omega  +ag_2\sigma\big)=0 \,.\label{lv2.3}
\end{eqnarray}
Taking the derivative with respect to $t$ in (\ref{lv2.3}) and using (\ref{lt2.3}) one obtains finally 
$$
\del^2_t \widetilde\omega -2a u_0\del^2_{x_1t}\widetilde\omega +u_0^2 \del_{x_1}^2\widetilde\omega +a g_2|D|\widetilde\omega =0
$$
or
\begin{equation}
(\del_t - a u_0\del_{x_1} )^2\widetilde\omega + u_0^2(1-a^2)\del_{x_1}^2 \widetilde\omega +ag_2|D|\widetilde\omega=0\,. \label{ellip}
\end{equation}
The nature of Equation (\ref{ellip}) depends on the parameters $a$, $g_2$ and $u_0$:
\begin{itemize}
\item If the Atwood number satisfies $a^2<1$ and if $u_0\neq 0$, then the equation
(\ref{ellip}) is always \emph{elliptic}, regardless of the sign of $ag_2$.
\item If $a^2=1$ and $ag_2>0$, which corresponds to the water waves case (with downwards gravity),  the equation (\ref{ellip}) is then \emph{nonstrictly hyperbolic} and the condition $ag_2>0$ is in fact the Levy condition ensuring the wellposedness of the Cauchy problem associated to it. 
\item If $u_0=0$ and $ag_2>0$, which corresponds to the Rayleigh-Taylor case  ($0<a<1$)
with the heavier fluid placed below the lighter one,  then (\ref{ellip}) has the same characteristics as in the case $a^2=1$, $ag_2>0$. The ill-posedness of the Rayleigh-Taylor
problem, even for configuration close to the rest state, are therefore due to \emph{nonlinear} effects.
\item If $a^2=1$ (or $u_0=0$) and $ag_2<0$, which corresponds to the water waves case (or Rayleigh-Taylor) with the lighter fluid \emph{below} the heavier one, the equation (\ref{ellip}) is then \emph{hypoelliptic}.
\end{itemize}
\begin{remark}
A Fourier type analysis of (\ref{ellip}) with the introduction of {\it modes}
\begin{eqnarray*}
\omega(x,t)= e^{\tau t +i\xi x}
\end{eqnarray*}
shows that all frequencies $\xi$ such that 
\begin{equation*}
|\xi|\ge \frac{a g_2}{(1-a^2) u_0^2}
\end{equation*}
lead  to unstable modes (with $\Re \tau>0$). The range of unstable frequencies increases with $a$ going to zero which is the Kelvin-Helmholtz case.  Without the nonlinear interactions, low frequencies would therefore remain stable for the Rayleigh-Taylor problem.
\end{remark}

\section{Consequences of the ellipticity of the Rayleigh-Taylor and Kelvin-Helmholtz problems}
\label{Csq}

\subsection{Analytic solutions}

For a precise statement of the behavior of analytic solutions to (\ref{nlgtangent})-(\ref{nlbs2}), it
is convenient to follow  Duchon and Robert  \cite{DR} and to introduce ${\mathcal B}_\rho$, the space of analytic functions which are Fourier transform of bounded measures with exponential decay at infinity (keeping in mind  the Payley Wiener theorem,  this corresponds to functions that can be analytically extended in the strip  $\{(x+ iy)\,, |y| <\rho\}$) :
\begin{equation}
{\mathcal B}_\rho= \{ u= \int d |\hat u|,\quad |\!| u|\!|_\rho =\int d |e^{\rho \xi}\hat u(\xi)|  <\infty\}.
\end{equation}
One then has the following theorem for the Kelvin-Helmholtz and Rayleigh-Taylor problems; it shows that
one can construct \emph{local} in time solutions to (\ref{nlgtangent})-(\ref{nlbs2}) if \emph{both} the interface
and the vorticity are prescribed (and analytic) at $t=0$, and that \emph{global} in time solutions exist
if the initial vorticity is not arbitrarily chosen. Note that in the second point of the theorem, one considers a system on $(\partial_{x_1}\sigma,\widetilde\omega)$ by taking the
$x_1$ derivative of (\ref{nlgtangent}).
\begin{theorem}
1- (see \cite{BFSS} and \cite{SS})  Let $\rho>0$ and  
$(\sigma_0, \widetilde \omega_0)\in {\mathcal B_\rho}$, and assume further that 
$\sigma_0$ has zero mean value (i.e.  $\widehat{\sigma_0}(\{0\})=0$).  Then
 there exists a  constant   $\alpha$ such that the Cauchy problem (\ref{nlgtangent})-(\ref{nlbs2})  has a unique  solution
 $$(\sigma(t,.),\tilde \omega(t,.))\in{\mathcal B}_{\frac {\alpha \rho-\abs{t}}{\alpha}}$$
 defined for $ |t| <\alpha \rho$ (and whose dependence on $t$ is analytic).\\
2-  (see \cite{DR}) There exists a constant $\epsilon>0$ such that for all 
  $\partial_{x_1}\sigma_0\in {\mathcal B}_0$ satisfying
\begin{equation}
   \int|\xi | d|\hat \sigma_0(\xi)| < \epsilon
\quad\mbox{ and }\quad
\widehat{\partial_{x_1}\sigma_0}(\{0\})=0,\label{duro}
\end{equation}
there exists $\alpha>0$ such that
the problem  $\partial_{x_1}$(\ref{nlgtangent})-(\ref{nlbs2})  admits  a unique solution  
$(\partial_{x_1}\sigma,\widetilde\omega)\in C(0,\infty ; {\mathcal B}_\alpha)^2 $ with the boundary conditions:
$$
\partial_{x_1}\sigma(0,x_1)=\partial_{x_1}\sigma_0(x_1),\qquad  \widehat{\widetilde{\omega}(0,\cdot)}(\{0\})=0.
$$
Moreover, $\partial_{x_1}\sigma(\cdot,t)$ and $\widetilde\omega(\cdot,t)$ converge to $0$ uniformly as $t$ goes to $\infty$.
\end{theorem}
\begin{remark}
\item[1-] The first point  has been proven first for the Kelvin-Helmholtz problem ($a=0$) in \cite{BFSS} and then extended to $a\not=0$ in \cite{SS}. The   basic new ingredient is contained in the relation between
$\tilde \omega$ and $\tilde \omega/2 - a(v_1+v_2\del_{x_1}\sigma)\,. $ The term $(v_1+v_2\del_{x_1}\sigma)$ is proportional to $v^\tau$ the tangent component of $v$ and one has:
$$
(v_1+v_2\del_{x_1}\sigma)=K(\tilde \omega)
$$
with $K$ being a regularizing operator as observed in Remark \ref{compacite}. Therefore  with Fredholm property and spectral analysis one can, for  $|a|<1$, invert the relation 
\begin{equation}
\frac{1}{2}\tilde \omega - a(v_1+v_2\del_{x_1}\sigma) =(\frac{1}{2}I-aK)(\tilde \omega)
\end{equation}
to obtain a well posed Cauchy problem for $(\sigma,\tilde \omega)$ in the class of analytic functions (see for details 
\cite{KL}, \cite{SS} and \cite{BHL}).
\item[2-] The proof of the second point
 has only been done for $a=0$ but the extension to $0\le|a|<1$ does not seem too difficult. 
\end{remark}
Looking at the solution provided by the second point of the theorem on the interval $[0,T]$, with $T>0$,
and changing the time variable $t\mapsto T-t$, we have at our disposal solutions which are analytic for 
$0\leq t<T$ and which at the time  $T$ exhibit any type of singularity compatible 
with the hypothesis  (\ref{duro}). 
This is a much larger class of singularities than the one initially predicted by Meiron, Baker and Orsag for 
Kelvin-Helmholtz \cite{MBO}. One can have examples for which the Hausdorff dimension of the graph 
$x_1\mapsto \partial_{x_1}\sigma(x_1,T)$ may be of dimension greater than one. In  particular in a periodic version of the above theorem one obtains starting, from analytic initial data at time $t=0$, 
a solution that coincides at  $t=T$ with the  Weierstrass function,
 $$
\sigma(x_1,T)=\sum_{0\le n \le \infty} 2^{-\frac n2 -n } \cos(2^n x_1),\qquad 
\widehat{\partial_{x_1}\sigma}(\xi,T)=\sum_{0\le n \le \infty} 2^{-\frac n2 } (\delta(\xi+ 2^n)+\delta(\xi-2^n)).
$$ 

The computations made in  \cite{MBO}  and \cite{M} led to the belief that the first appearance of singularity was the formation of a cusp in the curve   $\partial_{x_1}\sigma(x_1,T)$ or $\omega(x_1,T)$. 
As already said, this turns out to be a special case of the above theorem, but that was also the motivation of the construction of Caflisch and Orellana  \cite{CO} . With perturbation methods very close to the one used in   \cite{DR}, they construct solutions for Kelvin-Helmholtz problem which are analytic for $t>0$ and which at the time $t=0$ have an isolated singularity behaving like a cusp.  In complex representation these solutions are constructed as analytic perturbation of the function:
\begin{equation}
z_0(\lambda,t)=(1-i)\{(1-e^{-\frac t 2-i\lambda})^{1+\nu}-(1-e^{-\frac t 2+i\lambda})^{1+\nu}\},
\end{equation}
which turns out to be an exact solution of the linearized   Birkhoff-Rott equation and which has the following properties:
i) For any   $t>0$ $\lambda\mapsto z_0(t,\lambda)$ is analytic.
ii) For $t=0$ $\lambda\mapsto z(0,\lambda)$ does not belong to the space  $C^{1+\nu}$  but belongs to any space   $C^{1+\nu'}$ with $0<\nu'<\nu\,.$

\subsection{After the appearance of the first singularity}\label{after}

The mathematical properties of the interfaces may be used (with a contradiction argument) to obtain some results concerning the behavior of the solution (if it exists) after the appearance of the first singularities. As this turns out to be in many {\it free boundaries problems}, the method consists in showing that if the solution is locally (in space time) more regular than some threshold $\mathcal S$ then it is very regular (locally $C^\infty$ or analytic). This implies in turn that after the appearance of the first singularity, the solution has to be less regular than the threshold $\mathcal S$. As a consequence, the determination of a threshold $\mathcal S$ that would imply $C^\infty$ or analyticity is one of the challenges of the theory, and as such has been studied for instance by S. Wu \cite{Wu2} , 
 Lebeau \cite{LE} and Lebeau-Kamotski \cite{KL}. For instance in \cite{KL} one has the following statement: 
\begin{theorem}\label{KL}
 Assume that locally near a point   $(s_0,t_0)$ the curve $\Sigma(t)=\{r(s,t)\}$ parametrized by the arc length belongs to   $C^{1+\nu}$ ($\nu>0$) and is chord-arc; assume also that the density of vorticity is strictly positive and bounded at $(s_0,t_0)$. Then in the neighborhood of $(s_0,t_0)$, the curve and the density of vorticity are $C^\infty$. 
  \end{theorem}
\begin{remark}
If the vorticity is assumed to be positive and bounded for all $s$ (and not only in a neighborhood of $s_0$),
then the curve and density of vorticity are analytic. It is expected that the conclusion of Theorem \ref{KL} remains true if one replaces $C^\infty$-regularity by analyticity.
\end{remark}
 The result of  S. Wu \cite{Wu2} obtained independently and previously for Kelvin-Helmholtz  is slightly sharper. Instead  of $C^{1+\nu}\,$ one assumes that the function $r(t,s)$ belongs to $H^1$ and does not rolls up too fast. More  precisely one has the following theorem (we give a non sharp version of the condition (\ref{roll}) to avoid
the definition of technical functional spaces).
\begin{theorem} \label{wuth} Assume that ${\bf z}\in H^1([0,T],L^2_{loc}(\R))\cap L^2([0,T];H^1_{loc}(\R))$ 
is a solution  of the Birkhoff-Rott equation (\ref{BR}), which is uniformly chord-arc on $[0,T]$:
$$
\exists m,M>0,\quad \forall (t,\alpha,\beta)\in [0,T]\times\R^2,\qquad
m|\alpha -\beta| \le |{\bf z}(\alpha,t)-{\bf z}(\beta,t)|\le M|\alpha-\beta| .
$$
Then there exists a constant $c(m,M)$ such that the relation:
\begin{equation}
\sup_{[0,T] }|\!| \ln ({\bf z}_\alpha(\cdot,t))|\!|_\infty \le c(m,M)\label{roll}
\end{equation}
implies the analyticity.
\end{theorem}
  \begin{remark} The above statements are in full agreement with previous theoretical and numerical results.
  The example of   Caflish and Orellana corresponds to a singularity  in  $C^{1+\nu}$ and the Lebeau theorem shows that this function cannot be extended as a solution in  $C^{1+\nu'}$ with $\nu'>0$ after the singularity.\\
The  Prandlt-Munk vortex, described in Section \ref{weaksol} below, is a non regular solution of the Birkhoff-Rott equation defined for all time. This is not in contradiction with the previous theorems because near the points   $(\pm1,\frac t 2)$ the vorticity density does not satisfy the Lebeau-Wu hypothesis. Numerical simulations  done by Krasny \cite{KR} indicate the appearance after the singularity of spirals. On the other hand  Theorem \ref{wuth} implies that some spiral-type behavior would imply regularity and therefore is excluded. The catch lies in the fact that Krasny spirals do not seem to satisfy an hypothesis of the type (\ref{roll})
 \end{remark}

 Without entering into details it is interesting to give some comments on the proof of Theorem \ref{wuth}. As in any free boundary problem the challenge is to determine the threshold of regularity $\mathcal S$ 
as sharply as possible. The most delicate step is to prove that this threshold regularity implies a slightly better regularity; deducing analiticity (or $C^\infty$ regularity) from this slightly better regularity is
usually much easier. 
 
In the present case,  one shows that the assumptions of Theorem \ref{wuth} lead to $H^2$-regularity. This can be reached from the Birkhoff-Rott equation by a localized (with a smooth real valued function $\eta(\alpha)$) variational estimate. With  $V=\eta(\alpha)(1+i)\ln {\bf z}_\alpha$ one obtains the formulas
$$
\overline V_t-\frac{1}{2|{\bf z}_\alpha|^2}|D|V= F(\alpha,t),
$$
for some function $F(\alpha,t)$ that can be explicitly computed, and thus
\begin{equation}\label{Wu1}
\int  |{\bf z}_\alpha|^2| \big( \overline V_t-\frac{1}{2|{\bf z}_\alpha|^2}|D| V \big)^2d\alpha = \int |{\bf z}_\alpha|^2|F(\alpha,t)|^2 d\alpha.
\end{equation}
On the other hand, one easily checks that
\begin{eqnarray}
\int  |{\bf z}_\alpha|^2| \big( \overline V_t-\frac{1}{2|{\bf z}_\alpha|^2}|D| V \big)^2d\alpha&=&\int\big( \abs{{\bf z}_\alpha}^2\abs{V_t}^2+\frac{1}{4\abs{{\bf z}_\alpha}^2}\big\vert\abs{D}V\big\vert^2\big)\nonumber\\
\label{Wu2}
&-&\frac{1}{2}\Re \frac{d}{dt}\int V\abs{D}V.
\end{eqnarray}
Since the assumptions of Theorem \ref{wuth} imply that $m\leq \abs{{\bf z}_\alpha}\leq M$, it is possible to
deduce (after a bit of work!) from (\ref{Wu1})
and (\ref{Wu2}), a control of the quantity
$$
\int q(t)\int (m^2   |V_t|^2 + \frac 1{4 M^2}\big\vert|D|V\big\vert^2)d\alpha dt ,
$$
where $q(t)$ is a carefully chosen smooth, real valued function. Recalling that $V=\eta(\alpha)(1+i)\ln {\bf z}_\alpha$, this suggests that it is possible to control the $H^1_{loc}$-norm (in space and time) of $\ln {\bf z}_\alpha$ by
its $L^\infty$-norm, which is finite by assumption.

\section{Nonregular vortex sheets}\label{weaksol} 

\subsection{The Prandtl-Munk  vortex sheet}

As described in Section  \ref{weaksol0},  a solution of the Kelvin-Helmholtz problem (computed on the interface) should provide, with $u$ given by the Biot-Savart relation (\ref{biotsavart1}),  a weak solution of the incompressible Euler equation. The situation turns out not to be so simple for the following reasons:
\begin{enumerate}
\item Relaxing regularity  hypothesis on the interface, one can generate solutions of (\ref{prebs0})   which are not solutions of the Euler equation.  And this pathology appears with initial data which satisfy the hypothesis of \cite{LNZ} and therefore have weak solutions (different from what is given by the Kelvin-Helmholtz  equation (\ref{prebs0})).
\item Even without the above pathology, one should observe that the Kelvin-Helmholtz problem generates 
solutions which do not satisfy (and there is a big ``gap") the uniqueness stability criteria of Yudovich. 
\end{enumerate}
The point (i) is illustrated by the ``Prandtl-Munk vortex sheet" which  is notorious in fluid dynamics as the one that 
generates a circulation distribution that minimizes the induced drag on a plane wing 
\cite{Pu}. With initial data:
$$
 \Sigma(0)=\{r(x_1,0)=(x_1, 0)\},\qquad \widetilde\omega(x_1,0)=\frac{x_1}{\sqrt{1-x_1^2}}\chi_{(-1,1)}(x_1),
$$
where $\chi_{(-1,1)}$ is the characteristic function of the interval $(-1,1)$.
One observes   (cf. \cite{Pu}) that
  $$
 r(x_1,t)=(x_1, -\frac t 2),\qquad v  =(0,-\frac12),\qquad \widetilde \omega(x_1,t)= \widetilde\omega(x_1,0)
 $$
 is a solution of   the equation (\ref{prebs0}). 
 However the corresponding velocity  $u(x,y,t)$ is \emph{not} a solution of the Euler equation (\ref{Naples2})
and one has
instead (see \cite{LNS}):
  $$
 \nabla\cdot v=0,\qquad \del_t u + \nabla \cdot (u\otimes u)+\nabla p =F,
 $$
with a nonzero forcing term  $F$ given by (with $\delta$ denoting the Dirac mass):
 $$
  F(x_1,x_2)=\frac\pi 8  \Big((\delta_{(x_1 +1,x_2+\frac t2)}-\delta_{(x_1-1,x_2+\frac t2)}),0\Big).
$$
 
 {\it The catch lies in the fact that the vorticity density is not regular enough in the neighborhood of the points
 $x_1=\pm 1$ to justify the formal computations. On the other hand the initial data satisfy the hypothesis of \cite{LNZ} and there is therefore a weak solution of the Euler equations with these initial condition. This solution is singular enough near the same points to differ from the one given by  (\ref{prebs0}) .} We refer to
\cite{LLNS} for a discussion of the equivalence at low regularity of the Birkhoff-Rott and Euler description of vortex sheets.
 
\subsection{Regularization of vortex sheets}

There exist various ways to regularize a vortex sheet, with different properties.
\begin{enumerate}
\item Convolution by a sequence of functions converging to the Dirac distribution of  analytic initial data with vorticity concentrated on a curve  provides a regular solution of the Euler equation that does converge to the corresponding solution of the Kelvin-Helmholtz problem (in the existence domain of such solution), see \cite{MP}.
\item The so-called $\alpha$-regularization (one of the novel approaches for subgrid  scale modeling of turbulence): In $\R^2$ the vorticity transport equation:
is replaced by the $\alpha$ regularised  equation:
\begin{eqnarray} 
\del_t \omega_\alpha + u_\alpha \cdot \nabla \omega_\alpha  =0\,,  u_\alpha=K_\alpha \omega_\alpha= (I-\alpha \Delta )^{-1}u \,\hbox{ with } u(x,t)= \frac1{2\pi} R_{\frac \pi 2} \int_{\R^2} \frac{x-y}{|x-y|^2}\omega_\alpha(t,y)d y\,.\label{alphamodel}
\end{eqnarray}
For any $\alpha >0$ the kernel $K_\alpha$ is regularizing. Therefore the above system is well posed   (with global in time existence and uniqueness of the solution) for any measure valued initial data. Moreover any initial vorticity density supported by a curve $\Sigma(0)$ remains on the curve $\Sigma_\alpha (t)$ transported 
by the flow
$$
\dot X(x,t)=u_\alpha (X(x,t),t)\,,\,\, X(x,t=0)=x\,.
$$
Eventually if the initial data (curve and vorticity density) are analytic, for $\alpha \rightarrow 0$ the curve $\Sigma_\alpha(t)$ and the corresponding vorticity density converge to the solution of the Kelvin Helmholtz problem  as long as this solution, which exists locally in time, is defined; for details and further references see  
  \cite{BLT}.
 
\item However and more annoying is the following fact: with an initial vorticity concentrated on a curve the Navier-Stokes has with viscosity $\nu>0$ a well defined solution $u_\nu(x,t)$. Furthermore if the initial data satisfy the hypothesis of \cite{LNZ}   this solution converges to a weak solution. \\
 Weather this solution coincides with the one given by the Kelvin-Helmholtz problem (even with convenient hypothesis on the regularity of the curve and of the density of vorticity) is a completely open issue. 
 With the  ``Prandtl-Munk vortex sheet" one can exhibit a counter-example  (but with singular curve and vorticity) to this convergence:
 The corresponding initial data satisfy the hypothesis of \cite{LNZ}. Therefore, when the viscosity goes to $0$, the corresponding solution of the Navier-Stokes converges to a weak solution of the Euler equation. On the other hand the  ``Prandtl-Munk vortex sheet" which is the solution of the Kelvin Helmholtz equation with the same data is not a weak solution of the Euler equation.
 \end{enumerate}
 
 \section{Water waves}\label{swtt}
 
 The water waves problems refers to the case (like water in air) where one of the density (say $\rho_+$) is infinitely small with respect to the other and then taken equal to zero. 
As said in \S \ref{jump}, various parametrizations are possible. 
Parametrizing the surface $\Sigma=\{r(t,\lambda)\}$ as in (\ref{paramWW}), one has
 \begin{equation*}
a=-1,\qquad r_t= u_-^n \vec n+ u_-^\tau\vec \tau , \qquad r_\lambda={\vec \tau}|{r_\lambda}| \,.
 \end{equation*}
Therefore in  the system (\ref{tangent})-(\ref{vortex}), or equivalently
 \begin{eqnarray*}
&&(\del_tr -v)\cdot \vec n =0 ,\\
&& \del_t( \frac{\tilde\omega}2 + (v\cdot r_\lambda))+{\del_\lambda}\Big\{ \frac{1}{|r_\lambda|^2}(v-r_t)\cdot r_\lambda\,\big( \frac{\tilde\omega}2 + (v\cdot r_\lambda)\big)\Big\} 
+{\del_\lambda}\Big\{\frac{\tilde \omega^2}{8|r_\lambda|^2}-\frac{|v|^2} 2\Big\} -\vec {\mathbf g} \cdot r_\lambda= 0, 
\end{eqnarray*}
one may use the following relations
\begin{eqnarray*}
&&\frac{\tilde\omega}2 + (v\cdot r_\lambda)= u_-^\tau |r_\lambda |\,, \\
&&\frac{1}{|r_\lambda|^2}(v-r_t)\cdot r_\lambda\,( \frac{\tilde\omega}2 + (v\cdot r_\lambda))=(\frac12 u_-^\tau u_+^\tau-\frac{|u_-^\tau|^2}2) \,,\\
&&\Big\{\frac{\tilde \omega^2}{8|r_\lambda|^2}-\frac{|v|^2} 2\Big\}=(-\frac12 u_-^\tau u_+^\tau-\frac{|u_-^n|^2}2),
\end{eqnarray*}
to obtain
the equation
 \begin{equation}
 {\del_t}(u_-^\tau |{\del_\lambda r}|) -{\del_\lambda}\big (\frac 12(|u^n|^2+|u_-^\tau|^2 +\vec {\mathbf g}\cdot  r\big)=0\label{precurrent}\,.
 \end{equation}
A variant of (\ref{precurrent}) can also be derived if one chooses to parametrize the surface by its arclength (see for instance \cite{KL}); it is also worth commenting on the case
where the surface is a graph, $\Sigma(t)=\{(x_1,\sigma(x_1,t))\}$. Starting from 
(\ref{nlgtangent})-(\ref{nlgvortex}) and remarking that
$$
\frac{1}{2}\widetilde \omega+(v_1+v_2\partial_{x_1}\sigma)=\sqrt{1+(\partial_{x_1}\sigma)^2}u_-^\tau,
$$
one gets
$$
\partial_t \big(\sqrt{1+(\partial_{x_1}\sigma)^2}u_-^\tau\big)+\partial_{x_1}\Big\{u_-^\tau(v-r_t)\cdot\vec\tau\Big\}+\partial_{x_1}\Big\{\frac{1}{8}(u_+^\tau-u_-^\tau)^2-\frac{1}{2}\abs{v}^2\Big\}-\vec{\bf g}\cdot r_{x_1}=0.
$$
Since $r_t=(0,\partial_t\sigma)^T=(0,\sqrt{1+(\partial_{x_1}\sigma)^2}u_-^n)^T$, this is
equivalent to
\begin{equation}\label{etapB}
\partial_t \big(\sqrt{1+(\partial_{x_1}\sigma)^2}u_-^\tau\big)+\partial_{x_1}\Big\{
\frac{1}{2}(u_-^\tau)^2-u_-^\tau u_-^n-\frac{1}{2}(u_-^n)^2+g\sigma\Big\}=0,
\end{equation}
where we assumed that $\vec{\bf g}=(0,-g)^T$.\\
Since $\nabla\cdot u=0$ and $\nabla\wedge u=0$ in $\Omega_-(t)$, and assuming that 
$\Omega_-(t)$ is simply connected, there exists a velocity potential $\Phi(x,t)$ which 
 satisfies the relation:
 \begin{equation*}
 \nabla \Phi(x,t)= u(x,t),\qquad \Delta \Phi(x,t)=0
\quad\mbox{ in }\quad \Omega_-(t). 
\end{equation*} 
Denoting by $\Psi$ the trace of $\Phi$ at the surface, $\Psi(x_1,t)=\Phi(x_1,\sigma(x_1,t),t)$, one has
\begin{equation*}
\partial_{x_1}\Psi=\partial_{x_1}\Phi_{\vert_{\Sigma(t)}}+\partial_{x_1}\sigma\, \partial_{x_2}
\Phi_{\vert_{\Sigma(t)}} =\sqrt{1+(\partial_{x_1}\sigma)^2}\,u_-^\tau. 
\end{equation*}
Integrating with respect to $x_1$ in (\ref{etapB}) yields therefore
\begin{equation}\label{etapC}
\partial_t \Psi+
\frac{1}{2}(u_-^\tau)^2-u_-^\tau u_-^n-\frac{1}{2}(u_-^n)^2+g\sigma=0
\end{equation}
(note that adding if necessary a function depending on $t$ only to $\Phi$, the integration constant can be taken equal to $0$).\\
Eventually it is convenient to introduce the 
 Dirichlet Neumann operator $G(\sigma): H^{\frac12}(\R)\mapsto H^{-\frac12}(\R)$ associated to the resolution of the Dirichlet problem in $\Omega_- (t)\subset \Omega $ (with an additional homogeneous Neumann condition corresponding to (\ref{condN}) if there is a fixed boundary):
\begin{equation} G(\sigma)\Psi= \sqrt{1+(\partial_{x_1}\sigma)^2}\del_{\vec n}\Phi,  \label{neuman}
\end{equation}
where
$$
 \left\lbrace
\begin{array}{l}
-\Delta \Phi=0 \hbox{ in } \Omega (t),\\
 \Phi = \Psi \hbox{ on } \del \Omega_-(t)\backslash (\del\Omega_-(t)\cap \del\Omega ),\qquad
\del_{\vec n} \Phi =0  \hbox{ on }  (\del\Omega_-(t)\cap \del\Omega) .
\end{array}\right.
$$
Remarking that
$$
 \sqrt{1+(\partial_{x_1}\sigma)^2}u_-^\tau=\partial_{x_1}\Psi
\quad \mbox{ and }\quad  \sqrt{1+(\partial_{x_1}\sigma)^2}u_-^n=G(\sigma)\Psi,
$$
we deduce from (\ref{etapC}) that the water waves equations can we written
 \begin{eqnarray}
 \left\lbrace
\begin{array}{l}
\displaystyle \del_t \sigma=G(\sigma)\Psi,\\
 \displaystyle\del_t \Psi +g\sigma+\frac12(\del_{x_1}\Psi)^2-\frac1{2(1+(\partial_{x_1}\sigma)^2)}(G(\sigma)\Psi + \partial_{x_1}\sigma \partial_{x_1} \Psi)^2=0,
\end{array}\right.
\label{surfacewaves}
 \end{eqnarray}
where we assumed that $\vec {\bf g}=(0,-g)$ in the above expression.
\begin{remark} 
The formulation (\ref{surfacewaves}) of the water waves equation corresponds to the Hamiltonian exhibited by
Zakharov \cite{Zakharov} and was written first under this form by Craig and Sulem \cite{CrSu}. As it is well known the equation (\ref{surfacewaves}) can be obtained directly from the Euler equation
\begin{equation}
\del_t u + u\cdot \nabla u +\frac{1}{\rho_-}\nabla p =\vec {\bf g} \label{encoreueler}
\end{equation}
 written on the domain $\Omega_-(t)$ occupied by the fluid. The above short derivation has been written   to compare this problem with the Rayleigh-Taylor and Kelvin-Helmholtz equations and to emphasize similarities and differences.
\end{remark}

\subsection{Rayleigh-Taylor type behavior}

In the water waves case, the linearized equation (\ref{ellip}), with $u_0=0$, can be written
\begin{equation}\label{linww}
\del_t^2\widetilde\omega -g_2\abs{D}\widetilde\omega=0.
\end{equation}

It has already been noted in \S \ref{linan} that if $g_2>0$ (the fluid is above the surface), 
the operator appearing in (\ref{linww}) is hypoelliptic. The observations of Section \ref{Csq}
involving the ellipticity of the Kelvin Helmholtz problem could be adapted. {\it To the best of our knowledge this has not yet been worked out with full details except in the case of traveling waves which is the object of \S \ref{ertw}, and \S \ref{ertw2} for the extension to the $3d$ case}.
Near a stationary solution with $g_2>0$   one may prove, with analytic initial data, local in time existence of an analytic solution. Relaxing one  of the constraints at $t=0$ one may adapt the construction of \cite{DR} or \cite{CO}. Hence by  the change of variable $t\mapsto T-t$, one could construct solutions that  exhibit  a large class of singular behaviors. Finally since a Galilean transform changes the problem into an equation involving perturbations of a linear term of the form (\ref{linww}), the results of \cite{LE} and \cite{Wu2} could be adapted leading to a threshold of regularity $\mathcal S$ which for solutions more regular than this threshold would have $C^\infty$ or ``Gevrey" regularity. {\it In some sense the water waves problem with $ag_2<0\,$shares the same behavior as the Rayleigh-Taylor problem and therefore carries the same name.}

\subsection{The nonstrictly hyperbolic behavior}

When $g_2<0$ (the fluid is below the surface), we have already commented in \S \ref{linan} that (\ref{linww})
is nonstrictly hyperbolic, and that the condition $g_2<0$ corresponds to the Levy condition 
on the subprincipal symbol. The equation (\ref{linww}) roughly behaves as a wave equation with the positive self-adjoint operator $-\del_x^2$ replaced by its square root.

An important observation is that the sign of $g_2$ is an {\it intrinsic object}. The equation (\ref{linww}) has been obtained by small perturbation analysis near a point where the $x_1$ axis is tangent to the interface and where the interface behaves like $(x_1,\sigma(x_1,t))$ therefore $g_2$ coincides (modulo small perturbations) with $\vec {\bf g} \cdot \vec n$. With the full, nonlinear equations (\ref{encoreueler}), the corresponding quantity is $-\del_{\vec n} p$.
Therefore the condition 
\begin{equation}\label{RTstab}
-\del_{\vec n} p>0 \quad\mbox{ on }\quad \Sigma(t)
\end{equation}
should enforce the stability of the problem.  This condition is called the \emph{Rayleigh-Taylor stability criterion} (or \emph{Taylor sign condition}). From a mathematical viewpoint, it is also the \emph{Levy condition} ensuring the wellposedness of a nonstrictly hyperbolic equation. This is transparent on the formulation of the
water waves problem given in Proposition \ref{GMS1} below. Following \cite{GMS}, we denote by   $D_t= \del_t +v\cdot \nabla $ the lagrangian derivative, and by ${\mathcal N}$ the operator ${\mathcal N}=\sqrt{1+(\partial_{x_1}\sigma)^2}G(\sigma)$ (note that it is easy to define ${\mathcal N}$ even if $\Sigma(t)$ is not a graph); one then proves the following result for the interface curvature $\kappa$:  
\begin{proposition} \label{GMS1} (cf Proposition 2.6 in \cite{GMS} ).  For any smooth interface solution of the water wave equation one has:
\begin{equation}
D_t^2 \kappa +(-\del_{\vec n } p){\mathcal N} \kappa =R,  \label{courbure}
\end{equation}
where  $R$
is a remainder term consisting of lower order terms.
\end{proposition}
To prove this proposition one considers
\begin{equation}
\del_t u_- + u_-\cdot \nabla u_- +\frac{1}{\rho_-}\nabla p_- =\vec {\bf g} \label{encoreueler2}
\end{equation}
on the interface $\del\Omega_-(t)$
and apply  the Laplace Beltrami operator $\Delta_{\Sigma(t)}$ to  (\ref{encoreueler2}); using the important remark that at leading order, the Dirichlet-Neumann operator $G(\sigma)$ coincides with the square root of the Laplace Beltrami,
$$
G(\sigma)=\sqrt{-\Delta_{\Sigma(t)}}+R_0,\qquad
R_0=\mbox{zero-th order operator},
$$
one can establish that modulo lower order terms one has:
$$
\vec{n}\cdot \Delta_{\Sigma(t)}(\nabla p(t,x(t,\lambda))= \del_{\vec n} p\,{\mathcal  N} \kappa,
$$
and the result follows quite easily.

\medbreak

With the formula (\ref{courbure}) one obtains local in time stability of the interface for initial data having a limited Sobolev regularity. Of course, this result can be obtained with other techniques (see for instance \cite{Lindblad,WUb,La,CoutandShkoller,ShatahZeng}). 

\begin{remark}\label{remstr}
In the whole space $\R$  the solution of the model dispersive equation (corresponding to (\ref{linww}) with $g_2=-1$)
$$
\del_t^2u +  \abs{D}u=0,\qquad  u(x,0)=0,\qquad u_t(x,0)=\phi(x),  
$$
satisfies a dispersive estimate
\begin{equation}
\Vert u(\cdot,t)\Vert_{L^\infty}\le C t^{-\frac 1 2}\Vert \phi(x) \Vert_{L^1}. \label{strich}
\end{equation}
Therefore, one may expect to use dispersive effects, when the interface is a graph, to obtain {\it large time estimates} for the solution. 
This has been done in $2d$ (i.e. surface dimension $d-1=1$) by S. Wu \cite{Wu3}. The central idea in  \cite{Wu3} is to transform the equations (inspired by Birkhoff normal forms) in order to increase the degree in the non linear terms and therefore to reduce their size. Though considerable difficulties
make the implementation of this programme very delicate, it leads to an existence time of the order of 
$$
\exp \Bigg(\frac{C}{\hbox{some norm of the initial data}}\Bigg).
$$
Since the dispersive effect in the formula (\ref{strich}) improves with the dimension, the situation for this problem is better in $3d$ (i.e. surface dimension $d-1=2$, see \S \ref {ertw2}).
\end{remark}
\subsection{The Rayleigh-Taylor stability criterion}\label{rtcond}

Our goal here is to comment on whether the Rayleigh-Taylor stability criterion (\ref{RTstab}) is satisfied
for the water waves problem.
Let us consider  the case of water waves in finite depth. The boundary of $\Omega_-(t)$ is 
the union of two distinct pieces. A given boundary  $\Gamma_{fix}$ where the fluid is in contact 
with the boundary of the vessel, and the free boundary $\Sigma(t)$ between the region of the fluid  
of density $\rho_1$ (we take $\rho_-=1$ for the sake of simplicity),
and the fluid of density $\rho_+=0$. On $\Gamma_{fix}$ one assumes as usual the impermeability condition.
One thus gets
\begin{eqnarray}
&&\del_t u + \nabla (u\otimes u) = -\nabla p +\vec{\mathbf g},\qquad  \nabla \cdot u =0,\qquad
  \nabla \wedge u =0\quad \hbox{ in }\quad  \Omega_-(t),\label{int1}\\
 && p=0 \quad\hbox{ on }\quad \Sigma(t),\qquad
 u\cdot\vec n=0 \quad \hbox{ on }\quad   \Gamma_{fix}. \label{int2}
 \end{eqnarray}
Note that neither $\Sigma(t)$ nor $\Gamma_{fix}$ are assumed to be ``globally " a graph. 
However one consider only configurations where 
\begin{equation}
\label{nonvanishing}
d(\Gamma_{fix}, \Sigma(t))>h_{min}>0,
\end{equation} 
and therefore where the 
influence of the fixed boundary on the free boundary involves only lower order terms (in terms of regularity).

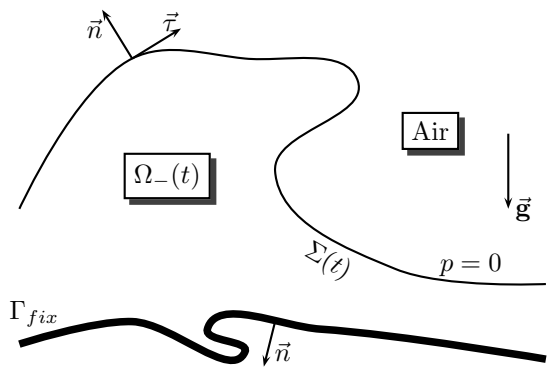
\begin{center}
\begin{pspicture}(7,5)
\pscurve[showpoints=false](0,2)(1.5,4)(3,4)(4.5,3.7)(3.4,2.5)(5,1.2)(7,1)
\psline{->}(1.5,4)(1.1,4.65)
\rput(1,4.4){$\vec n$}
\psline{->}(1.5,4)(2.15,4.4)
\rput(2,4.5){$\vec\tau$}
\rput(5.5,3){\psshadowbox{Air}}
\rput(2,2.4){\psshadowbox{$\Omega_-(t)$}}
\rput{330}(4.1,1.3){$\Sigma(t)$}
\rput(6,1.25){$p=0$}
\pscurve[showpoints=false,linewidth=.1](0,0.2)(1.5,0.5)(3,0.1)(2.5,0.4)(4,0.4)(7,0)
\rput(0.2,0.6){$\Gamma_{fix}$}
\psline{->}(6.5,3)(6.5,2)
\rput(6.7,2){$\vec {\bf g}$}
\psline{->}(3.4,0.5)(3.25,-0.1)
\rput(3.5,0.1){$\vec n$}
\end{pspicture}
\end{center}

 Taking in particular into account the relation $\nabla\wedge u_-=0$, 
one deduces from (\ref{int1})-(\ref{int2}) that the pressure $p$ satisfies 
$$
\left\lbrace
\begin{array}{l}
-\Delta p= \sum_{ij}\del_{x_i}u_j\del_{x_j} u_i = |\nabla u|^2\ge 0 
\quad \hbox{ in }\quad \Omega_-(t),\vspace{1mm}\\
 \displaystyle p=0 \quad \hbox{ on }\quad  \Sigma(t),\qquad \del_{\vec n} p= -\sum_{ij} u_i u_j \frac{\del {\vec n_i}}{\del x_j} +\vec{\mathbf g}\cdot \vec n \quad \hbox{ on }\quad  \Gamma_{fix}.
\end{array}\right.
$$
Then the maximum principle implies that   $p> 0$ in  $\Omega_-(t)$ and therefore that
$$
-\del_{\vec n} p> 0 \quad \hbox{ on }\quad\Sigma(t),
$$
as soon as
\begin{equation}
\del_{\vec n} p= -\sum_{ij} u_iu_j \frac{\del {\vec n_i}}{\del x_j} +\vec{\mathbf g}.\vec n> 0  \quad \hbox{ on }\quad  \Gamma_{fix}.\label{pstab}
\end{equation}
For instance (\ref{pstab})  turns out to be true for flat bottoms, since one then has
$$
 \sum_{ij} u_iu_j \frac{\del {\vec n_i}}{\del x_j} +\vec{\mathbf g}\cdot\vec n=\vec{\mathbf g}\cdot\vec n> 0;
$$
the same argument works also in infinite depth (in fact, this is the configuration where it has been proved first \cite{WU,WUb}), or for liquid drops\footnote{However, this argument fails if the fluid is not irrotational. This is the reason why, in the study of liquid drops, it is \emph{assumed} in \cite{Lindblad,ShatahZeng,CoutandShkoller} that the Rayleigh-Taylor condition holds. It also known that without this condition, the equations are ill-posed \cite{Ebin}}.

\begin{remark}
Flat, infinite bottoms, or liquid drops are of course not the only configurations where the Rayleigh-Taylor stability criterion (\ref{RTstab}) is satisfied. Quite obviously, the same arguments show that (\ref{RTstab}) is also
satisfied for ``almost'' flat bottoms. More generally, a condition on the curvature of the bottom can
be derived   \cite{La,AL}.
 \end{remark}

\subsection{Ellipticity and regularity of traveling waves}\label{ertw}

Ellipticity results have been used above to derive instability results for the 
Kelvin-Helmholtz and Rayleigh-Taylor problems. Though such instabilities do not occur 
for surface waves (when the Taylor sign condition (\ref{RTstab}) is satisfied), an elliptic behavior can still be observed in some phenomena concerning water-waves. In particular, it is possible to derive \emph{a priori} regularity results for traveling waves, both in $2d$ and in $3d$. Comments on the $3d$ case are postponed
to \S \ref{ertw2}.

Using the formulation (\ref{surfacewaves}) of the water-waves in terms of the surface elevation $\sigma$ and the trace of the velocity potential at the free surface $\Psi$, a traveling
wave is a solution to (\ref{surfacewaves}) of the form
$$
(\sigma(t,x_1),\psi(t,x_1))=(\underline{\sigma}(x_1-c t),\underline{\psi}(x_1- ct)),
$$ 
where $c\in \R$ is the speed of the traveling wave. 
Consequently,
$(\underline{\sigma},\underline{\psi})$ must solve (we omit the underlines for the sake of clarity),
\begin{equation}\label{formA}
\left\lbrace
\begin{array}{l}
\displaystyle G(\sigma)\psi+c\partial_{x_1} \sigma=0,\\
\displaystyle g\sigma-c\partial_{x_1}\Psi+\frac{1}{2}\vert\partial_{x_1}\Psi\vert^2-\frac{1}{2}\frac{(\partial_{x_1}\sigma\partial_{x_1}\Psi-c\partial_{x_1}\sigma)^2}{1+\vert\partial_{x_1}\sigma\vert^2}=0.
\end{array}\right.
\end{equation}

Existence of traveling waves in one surface dimension is known from the early 
works of Levi-Civita \cite{LeviCivita} and Nekrasov \cite{Nekrasov}. The fact that if the free boundary of these solutions is a $C^1$ curve, then it is also analytic is also a classical result that goes back to the works of Lewy \cite{Lewy} and Gerber \cite{Gerber1,Gerber2,Gerber3}. These results use the fact that the system (\ref{formA}) is
\emph{elliptic}, a fact that can easily be observed on the linearized equations around the rest state $\sigma=0$, $\psi=0$. Recalling that $G[0]=\vert D \vert$, the
linearized equations read indeed
$$
\left\lbrace
\begin{array}{l}
\displaystyle \abs{D}\psi +c\partial_{x_1} \sigma=0,\\
g\sigma-c\partial_{x_1}\psi=0.
\end{array}\right.
$$
The symbol of this system is therefore 
$$
\left(\begin{array}{cc}
  ic\xi_1 & \abs{\xi_1}\\
g & -ic\xi_1
\end{array}\right),
$$
whose determinant is the elliptic symbol $c^2\xi_1^2-g\abs{\xi_1}$.

\subsection{A word on shallow water asymptotics}\label{shallow}

Except for its influence, through (\ref{pstab}), on the Rayleigh-Taylor criterion (\ref{RTstab}), the bottom
does not play any role in the well-posedness theory of the water waves problem described above. By assuming
(\ref{nonvanishing}), i.e. that the water depth never vanishes, one is ensured that the contribution
of the bottom to the Dirichlet-Neumann operator 
$G(\sigma)\Psi$
is \emph{analytic} (by the ellipticity of the potential equation $\Delta\Phi=0$); therefore, \emph{the contribution of the bottom is ignored in any symbolic description of $G(\sigma)\Psi$}. This is the reason why it is possible to obtain a local well-posedness theory for very rough bottoms, as shown recently in \cite{ABZ}.

For the analysis of instabilities that interests us here, the influence of the bottom is thus irrelevant. 
This fact seems to contradict common sense observations, which show that the bottom often plays a important role in the dynamics of the surface (consider the shoaling of a wave as it approaches the shore for instance). The explanation of this apparent paradox is that the fact that the contribution of the bottom is \emph{smooth} (and can therefore be ignored in the wellposedness theory) does not mean that it is \emph{small}. For instance, consider
the approximation of the Dirichlet-Neumann operator by the square root of the Laplace Beltrami of the surface,
$$
G(\sigma)=\sqrt{-\Delta_{\Sigma(t)}}+R_0,\qquad
R_0=\mbox{zero-th order operator};
$$
we have already seen that this symbolic approximation is very useful to prove the stability of the water waves problem. However, it is not uniform in the shallow water limit. More precisely, the operator norm $\Vert R_0 \Vert_{L^2\to L^2}$ grows to infinity in the shallow water limit (which corresponds to configurations where the ratio of the depth over the typical wavelength of the interface deformation is small). In order to obtain uniform estimates in the shallow water limit, one must take into account 
the contribution of the topography\footnote{In \cite{LanIW}, a symbolic analysis ``with tail'' of the Dirichlet-Neumann operator is proposed, the tail being an infinitely smoothing perturbation of the principal symbol $\sqrt{-\Delta_{\Sigma(t)}}$ that takes into account the
contribution of the bottom and removes the shallow water singularity}.   In some sense, the study of the interface is not decorrelated from the topography in the shallow water regime (the same phenomenon also occurs for interface waves, see for instance \cite{BLS}). In order to focus on the instability phenomena at the interface, 
\emph{we implicitly assumed throughout this paper that we are \emph{not} in a shallow water type regime}. The study of this latter is however possible, albeit with different techniques, and a consequent literature has addressed this problem. A few examples are \cite{Craig,KN,SW,Iguchi,AL}.
 
 \section{ Differences between the $2d$ the $3d$ case}\label{3d}

 Most of this review has been devoted to interfaces in planar flows. Before concluding and without going to the details one should single out differences and similarities between the $2d$ and the $3d$ case.
 
 In $3d$ the  formula for the propagation of the vorticity  is 
 $$
  \del_t \omega +u\cdot \nabla \omega = \omega \cdot \nabla u;
 $$
 as a consequence the $2d$ conservation relation does not exist any more. 
 Therefore, there is no theorem for weak solutions nor global in time result for strong  solutions. The only general available result concerns local existence in time  of smooth solution, say,  in $C^{1,\alpha}$.  This
regularity seems optimal as shown by the instability of the shear flow \cite{Majda-Bertozzi,BT}:
more precisely,  consider the vector field
 \begin{equation}
 U_{\rm shf}(x,t)=(u_1(x_2), 0, u_3(x_1-tu_1(x_2)),
 \end{equation}
 which is an obvious solution of the $3d$  incompressible Euler equation. With this type of flow one can construct solutions with the following property
 $$
 U_{\rm shf}(x,0)\in C^{0,\sigma}, \quad \hbox{ and for any } \quad t>0,\qquad
 U_{\rm shf}(x,t)\notin C^{0,\sigma}\,.
$$
\subsection{ About the $3d$ Kelvin-Helmholtz and Rayleigh-Taylor problems}

For smooth solutions, the generalization of the formulas (\ref{tangent})-(\ref{vortex}) for the equations of the interface remain valid. Hence one can prove a local in time existence theorem for analytic initial data (cf. \cite{SS}). 
 
However the ellipticity properties cease to be valid. As previously observed, these properties follow from the fact that the equation can be locally considered as a perturbation of the linearized problem. This 
linearized problem  is elliptic in $2d$. In $3d$, one can still 
reduce the problem (see \cite{BFSS,CHA})
to the situation where the interface is  a graph  $x_3=\sigma(x_1,x_2,t)$,
  and where the solutions are perturbations of the stationary state  $\underline{\sigma}=0$, $\widetilde{\underline{\omega}} (x_1,x_2)= ( \widetilde{\underline{\omega}}_1, \widetilde{\underline{\omega}}_2,0)$.
   Then in Fourier variables the linearized problem in the Kelvin-Helmholtz case (the Rayleigh-Taylor case can be studied in the same way) is
$$
 \del_t \left(\begin{array}{c} \widehat {\sigma}\\   \widehat {  \omega}_1 \\ \widehat{  \omega}_2\\ \widehat{  \omega}_3\end{array}\right)={\mathcal A}(\xi)\left(\begin{array}{c} \widehat {\sigma}\\   \widehat {  \omega}_1 \\ \widehat{  \omega}_2\\ \widehat{  \omega}_3\end{array}\right),
 $$
where the operator ${\mathcal A}(\xi)$ is given by (denoting $\xi=|\xi| (\cos \theta, \sin \theta)$) 
 \begin{equation}
 {\mathcal A}(\xi)=\left(\begin{array}{clcr} 0 &\frac i2 \sin \theta  & -\frac i2 \cos \theta  & 0 \\
- \frac i2 |\xi|^2|\widetilde{\underline{\omega}}|^2\sin \theta &0  &0 &\frac12(\xi\cdot\widetilde{\underline{\omega}})\sin \theta\\
\frac i2 |\xi|^2|\widetilde{\underline{\omega}}|^2\cos \theta &  0 &0 &-\frac12(\xi\cdot\widetilde{\underline{\omega}})\cos \theta\\
0 &-\frac12(\xi\cdot\widetilde{\underline{\omega}})\sin \theta  &\frac12(\xi\cdot\widetilde{\underline{\omega}})\cos \theta  & 0
\end{array}\right)\,.
 \end{equation} 
The eigenvalues of the matrix   ${\mathcal A}(\xi)$ are
$$0 \mbox{ (with multiplicity $2$), }\quad\mbox{ and }\quad \pm \frac12|\xi\wedge \widetilde{\underline{\omega}}|.$$
Therefore, the operator 
$\del_t -{\mathcal A}(\xi)$ is no more elliptic and perturbations transverse to the direction of streaming are thus unaffected \cite{CHA}.

\subsection{About $3d$ water waves}\label{ertw2}

\subsubsection{A global existence result in $3d$}

Local existence for $3d$ (surface dimension $d-1=2$) water waves under  
the Rayleigh-Taylor criterion (\ref{RTstab}) can be proved with tools that do not
depend on the dimension \cite{WUb,La,CoutandShkoller,ShatahZeng}. This is of course not the case if one
wants to prove a global existence result based on the dispersive properties of the equation.\\
In order to describe the stability result for the water waves
 in $2d$ (surface dimension $d=1$), we used the equation on the curvature $\kappa$
provided by Proposition \ref{GMS1}. In \cite{GMS}, it is proven that such an equation still holds in $3d$ with the curvature replaced by the mean curvature, and a careful analysis of this equation allows the authors
to obtain a global existence result for small datas.  This is achieved by using several ingredients
(in the same spirit as those in $2d$, see Remark \ref{remstr}, but with different techniques):
\begin{itemize}
\item In the whole space $\R^2$,  the solution of the model dispersive equation
$$
\del_t^2u + \abs{D}u=0, u(x,0)=0 \,, u_t(x,0)=\phi(x)\,,                        
$$
satisfies a decay estimate
\begin{equation}
\Vert u(x.,t)\Vert_{L^\infty}\le C t^{-1}\Vert\phi \Vert_{L^1}.
\end{equation}
\item The use of a Birkhoff normal form to increase the degree of the non linearity. However, time resonances
do not allow the construction of a standard normal form; the authors therefore introduce the notion of space resonance, which catches the fact that, for dispersive equations, wave packets do not have the same group velocity and therefore ``do not see each other'' in general.
\end{itemize}

\subsubsection{Traveling waves in $3d$}
Traveling waves also exist for $3d$ water waves (surface dimension $d-1=2$), but this fact is considerably more difficult to establish than in the $2d$ case and has only been solved very recently by Iooss and Plotnikov \cite{IoossPlotnikov}. As in the $2d$ case, a simple
look at the linearized equations is very instructive. The symbol of the linearized
system is indeed given by
$$
\left(\begin{array}{cc}
  ic\xi_1 & \abs{\xi}\\
g & -ic\xi_1
\end{array}\right),
$$
with $\abs{\xi}=(\vert \xi_1\vert^2+\vert\xi_2\vert^2)^{1/2}$. The determinant is then the symbol $c^2\xi_1^2-g\vert \xi\vert$. This symbol is \emph{not elliptic} since
it vanishes on the (unbounded) region ${\mathcal C}=\{c^2\xi_2^2=g\abs{\xi}\}$. The
equation defining ${\mathcal C}$ is a Schr\"odinger equation on the boundary that may propagate singularities.

 Despite of this, Alazard and M\'etivier \cite{AlazardMetivier} managed to prove that smooth enough biperiodic (i.e. periodic in both space directions) traveling waves are automatically $C^\infty$.
By a careful paralinearization of the equations and two transformations inspired by the work of Iooss and Plotnikov, they show that the general situation
is quite similar to what is observed in the linearization around the rest state.  Then, they exhibit a diophantine condition that allow them to handle the degeneracy of the principal symbol.

\begin{remark}
As already mentioned, the sign of gravity also plays an important role. Were it negative, the symbol of the linearized equation would be weakly elliptic, even in the $3d$ case. This remains true for the full equations (\ref{formA}) if the Rayleigh-Taylor condition is not satisfied. The paralinearization strategy of Alazard-M\'etivier can then be used.
\end{remark}

\section{Conclusion}
As a conclusion, we  recall some of the basic threads that lead to the organization of this survey and  add extra remarks and open problems.

All the analysis presented here is deduced from the incompressible Euler equation written in the sense of distributions. With this setting the stability of the Cauchy problem has been analyzed in the models problems 
of Rayleigh-Taylor (Atwood number $0<\abs{a}< 1$), Kelvin-Helmholtz ($a=0$), and water waves   ($\abs{a}=1$). 

With an high frequency analysis it has been shown that for $|a|<1$ the linearized problem is elliptic. Therefore the Cauchy problem is well posed only with analytic initial data. 
Moreover the solution whenever it exists is analytic even if the data (say for $t=0$ or $T=t$ ) exhibit a singular behavior. 
To persist after the singularity the interface has to be very {\it wild } and this may be in agreement both with the numerical simulation of \cite{KR} and with physical observations.

The situation for the water waves may be different. When the Rayleigh-Taylor sign condition $-\del_{\vec n} p>0$ is not satisfied the linearized problem shares many similarities with the Rayleigh-Taylor equations. 
On the other hand with $-\del_{\vec n} p>0$ on the interface stability results can be obtained. 

The above results given in $2d$ are compared in Section \ref{3d}  with the situation in $3d$. The ellipticity (or hypoellipticity) of the underlying linearized problem is lost and therefore the situation 
may be very different for the Kelvin-Helmholtz and Rayleigh-Taylor problems. 
As a companion problem observe that with the shear flow  
$$u(x_1,x_2,x_3,t)=(u_1(x_2), 0, u_3(x_1-tu_1(x_2)),$$
one can construct solutions defined for all time and with any ``limited" degree of regularity on the interface.
The same kind of difficulty is also encountered in some situation for the water waves problem (for the regularity of traveling waves for instance).

\medbreak

In the present description physical quantities like the surface tension and the viscosity have been ignored.  
This is based on the observation that such quantities are {\it small}\footnote{this assertion is not very honest! Viscosity and surface tension are negligible for physical phenomena involving ``reasonable'' frequencies, but most of the instabilities described here occur at very high frequencies; the viscosity and surface tension terms should then play an important role. We refer to \cite{LanIW}, where it is shown that a natural two-fluid generalization of the Rayleigh-Taylor stability criterion (\ref{RTstab}) is given by
$$
[-\partial_z P]>\frac{1}{4}\frac{(\rho^+\rho^-)^2}{T(\rho^++\rho^-)}{\mathfrak c}(\sigma)\vert\omega\vert_\infty^4,
$$
where $T$ is the surface tension coefficient, $\omega$ is the jump of the horizontal velocity at the interface, and ${\mathfrak c}(\sigma)$ is a geometrical constant.\\
However this ``not very honest mathematical behavior" is not restricted to fluid mechanics and nevertheless contributes to the understanding of nature. An other basic example is the Maxwell equations; they involve the permittivity $\epsilon(x)$ and  the magnetic permeability $\mu(x)$. These equations are the prototype of hyperbolic systems  and the high frequency analysis is systematically used for their solutions. On the other hand $\epsilon$ and $\mu$ are derived from the Helmholtz equation under the hypothesis that the frequency is ``moderate" . For higher frequencies $k$ in the Helmholtz equation they would depend on $k$. Therefore in the space domain one obtains an equation with $\epsilon$ and $\mu$ acting as time convolution and no more as coefficients and the nature of the equations at high frequency should be changed!
} and that the understanding of the 
pathology of the phenomena when they are made equal to $0$ should contribute to the understanding of 
what happens when they are small. 
\begin{itemize}
\item The introduction of a surface tension  creates a third order dispersive term and
all the above equation become locally in time well posed \cite{ITT,Am,Ammas}. This property does not persists when this surface tension goes to zero 
except  for the water waves problem with the Rayleigh-Taylor sign condition
\cite{AM2,ShatahZeng,MZ}.
\item With the presence of a viscosity $\nu>0$ the Euler equation is changed into  the Navier-Stokes equation. When the density is constant ($\rho_+=\rho_-$), this equation is in $2d$ reasonably
understood even for initial data with a vorticity being a measure concentrated on a curve. When $\nu $ goes to zero (under reasonable hypothesis on the sign of the vorticity $\nabla\wedge u_\nu$) the corresponding solution converges to a weak solution of the Euler equation. 
Since there is no uniqueness theorem for weak solutions of the Euler equation, the connection between this limit and corresponding solutions of the Kelvin Helmholtz equation remains a fully open problem. Moreover with the Prandlt-Munk vortex described in Section \ref{weaksol} one can show the existence of situations where such convergence does not hold.
\item 
Eventually one could consider a viscous perturbation of  the water waves problem, i.e. a Navier-Stokes equation with a free boundary. This has already been studied (cf. for instance \cite{Sol}). One may conjecture that in this case, with the Rayleigh-Taylor sign condition at the boundary, the convergence to the water waves solution should be true.
\end{itemize}

{\it Acknowledgments.}
Ce travail a b\'en\'efici\'e d'une aide de l'Agence Nationale de la Recherche portant  la r\'ef\'erence ANR-08-BLAN-0301-01

\end{document}